\documentclass[a4paper]{amsart}

\usepackage{amssymb}
\usepackage{amsthm}  
\usepackage{amsmath} 
\usepackage{amscd} 
\usepackage[all]{xypic}
\usepackage{url}

\newcommand{\omi}{\omega^{-1}}
\newcommand{\om}{\omega}
\newcommand{\ka}{\kappa}
\newcommand{\si}{\sigma}

\newcommand{\ba}{\mathcal{G}} 
 
\newcommand{\va}{\varphi}
\newcommand{\fg}{\mathfrak g}
\newcommand{\fp}{\mathfrak p}
\newcommand{\Ad}{{\rm Ad}}
\newcommand{\Fl}{{\rm Fl}}
\newcommand{\Aut}{{\rm Aut}}
\newcommand{\id}{{\rm id}}
\newcommand{\na}{\nabla}
\newcommand{\U}{\Upsilon}

\newcommand{\Rho}{{\mbox{\sf P}}}
\newcommand{\R}{\mathbb{R}}

\newcommand{\Om}{\Omega}

\newtheorem*{prop*}{Proposition}

\newtheorem*{thm*}{Theorem}

\newtheorem*{lemma*}{Lemma}

\newtheorem*{cor*}{Corollary}

\newtheorem*{rem*}{Remark}
\theoremstyle{definition}
\newtheorem*{def*}{Definition}

\theoremstyle{remark}
\newtheorem*{exam}{Example}

\begin{document}
\title{Parabolic Symmetric Spaces
}
\author{Lenka Zalabov\' a}
\address{
Eduard \v Cech Center for 
Algebra and Geometry, Masaryk University, 
Faculty of Science, Kotl\' a\v rsk\'a 2,
611 37 Brno, Czech Republic
}
\email{zalabova@math.muni.cz}
\dedicatory{This paper is dedicated to Peter Michor on the occasion of his 60th birthday.} 
\thanks{
This research has been supported at different times 
by the ESI Junior Research Fellowships program of The Erwin Schr\" odinger 
International Institute for 
Mathematical Physics and by the 
Eduard {\v{C}}ech Center for Algebra and Geometry, project nr.\ LC505. 
The author acknowledges very useful discussions with Andreas \v Cap 
during the work on this paper.
}
\begin{abstract}
We study here systems of symmetries on $|1|$--graded parabolic geometries. 
We are interested in smooth systems of symmetries and we discuss non--flat homogeneous
$|1|$--graded geometries. We show the existence of an invariant admissible 
affine connection under quite weak condition on the system.
\end{abstract}
\maketitle
Affine symmetric spaces are well known objects in differential geometry. 
One can find classical description of them for instance in \cite{H} or \cite{KN2}.
Let us sketch here briefly this classical concept:
Let $M$ be a manifold with an affine connection $\na$. A \emph{symmetry} at $x$ is a 
(globally defined) diffeomorphism $s_x$ of $M$ such that
$s_x(x)=x$, $T_xs_x=-\id$ on $T_xM$ and $s_x$ is an affine transformation of $\na$. 
Equivalently, a symmetry at $x$ is an involutive affine transformation with 
isolated fixed point $x$.
Clearly, there can exist at most one symmetry $s_x$ at each 
$x$ on a manifold $M$ with a connection $\na$ and a pair $(M,\na)$ is 
called \emph{symmetric space}, if there exists a symmetry at each $x \in M$. 
Each symmetric space is homogeneous because the group of affine transformations 
(involving symmetries) acts transitively on $M$. 
Let us point out here that on each symmetric space $(M,\na)$, 
uniqueness of symmetries implies that $s_x \circ s_y \circ s_x=s_{s_x(y)}$ 
holds for each $x,y \in M$ and in fact, we can describe the (uniquely 
given) system of symmetries as a smooth multiplication
$S:M \times M \rightarrow M$, $(x,y)\mapsto s_x(y)$.

Following \cite{L,K}, one can also define a symmetric space algebraically 
as the following structure: A \emph{symmetric space} is a manifold $M$ together
with a smooth multiplication $S:M \times M \rightarrow M$, $S(x,y)=s_x(y)$ 
with the following four properties: $s_x (x)=x$,  $s_x \circ s_x (y)=y$,    
$s_x \circ s_y(z)=s_{s_x(y)} \circ s_x(z)$ and every $x$ has a neighborhood 
such that $s_x(y)=y$ implies $y=x$ for each $y$ from the neighborhood.
Surprisingly, this definition is equivalent with the latter description. Really, 
in \cite{L} O. Loos proved that each symmetric space $(M, S)$ admits a unique affine 
connection $\na$ which is invariant with respect to all symmetries $s_x$. 
Hence the structure $(M,S)$ generates a unique affine symmetric space $(M, \nabla)$.
This point of view of symmetric spaces is the most interesting for us. 
Let us remark that, instead of the fourth property in the definition, we can suppose 
$T_xs_x= -\id$ on $T_xM$, which makes the
definition by O. Loos more compatible with our basic definitions, see \ref{2.1}.

In this article, we discuss systems of symmetries 
which instead of an affine connection preserve a $|1|$--graded 
parabolic geometry. These form a rich family of geometries which involves well known 
examples like conformal or projective structures. Following the classical point of view, 
we define a symmetry on a $|1|$--graded geometry as a morphism satisfying 
$s_x(x)=x$ and $T_xs_x=-\id$ 
on $T_xM$ and we require that the symmetry is a morphism 
of the $|1|$--graded geometry and thus preserves the corresponding geometric 
structure (e.g. conformal or projective). 
In this definition, we do not assume it is an affine transformation.

Contrary to the affine case, there can exist many different symmetries at one point 
on a $|1|$--graded geometry. It comes from the fact that parabolic geometries are 
structures of second order and our definition of the symmetry prescribes only 
$1$--jet of the morphism, see \cite{ja-dis, ja-elsevier}.
Consequently, there can exist various systems $S: M \times M \rightarrow M$, 
$S(x,y)=s_x(y)$ of symmetries on a $|1|$--graded geometry and such a 
system $S$ does not necessarily carry an invariant affine connection. 

In this article, we look for conditions under which $S$ is a system of 
affine transformations of a suitable admissible affine connection. 
The main motivation for us is the article \cite{P}, where the author discusses the 
projective case in the classical setup. The general theory of parabolic geometries 
allows us to discuss all $|1|$--graded geometries in a uniform way. 
More precisely, we use the theory of Weyl structures (see \cite{WS, parabook}) 
to find an invariant connection with properties as above under quite weak 
conditions on the system $S$.  We will see that the conditions from the Loos' algebraic 
definition of a symmetric space are crucial here.
Our first condition on the system $S$  is its smoothness. Contrary to the affine 
symmetric spaces, there can exist non--smooth systems of symmetries for general 
$|1|$--graded geometries. Assuming the smoothness of $S$, 
we concentrate on non--flat homogeneous $|1|$--graded geometries.
The second and last condition we need is exactly 
the condition $s_x \circ s_y \circ s_x=s_{s_x(y)}$ from the 
Loos' definition. We show that if the two conditions on $S$ are satisfied, then the 
$|1|$--graded geometry reduces 
to an affine symmetric space (for the system $S$). We also show that in many cases, 
these conditions are trivially satisfied and then there are no other 
examples than the affine symmetric spaces.

\section{$|1|$--graded parabolic geometries}
In this introductory section, we remind roughly basic definitions 
and facts on Cartan and parabolic geometries.
In this paper, we follow concepts and notation of \cite{WS, parabook} 
and the reader can find all details therein. We are mainly interested in 
$|1|$--graded parabolic geometries. The complete list of them can be found
in \cite{WS, ja-dis}.

\subsection{Parabolic geometries} \label{1.1}
Let $G$ be a Lie group, $P\subset G$ its Lie subgroup, and denote by 
$\fp\subset\fg$ corresponding Lie algebras. 
A \emph{Cartan geometry} of type $(G,P)$ on a smooth manifold $M$ 
is a pair $(\ba\rightarrow M,\om)$ consisting of a
principal  $P$--bundle $\ba \rightarrow M$ and of a
one--form $\om \in \Omega^1(\ba,\fg)$, called the \emph{Cartan connection}, 
which is $P$--equivariant, reproduces generators of fundamental 
vector fields and induces a linear isomorphism $T_u\ba\cong\fg$ for each $u \in \ba$. 
The $P$--bundle $G\rightarrow G/P$ together with the (left) Maurer--Cartan form
$\om_G \in \Om^1(G,\fg)$ forms a Cartan geometry of type $(G,P)$ 
which is called \emph{homogeneous model} or \emph{flat model}. 

A \emph{morphism} between Cartan geometries of type $(G,P)$ from 
$(\ba \rightarrow M,\om)$ to $(\ba'\rightarrow M',\om')$ is a $P$--bundle 
morphism $\va:\ba \rightarrow \ba'$ such that $\va^*\om'=\om$. 
Further we denote the base morphism of $\va$ by $\underline \va:M \rightarrow M'$.
For simplicity, we suppose that the maximal normal subgroup $K$ of $G$ 
which is contained in $P$ is trivial. 
With this assumption, there is one to one correspondence between the morphisms 
$\va$ and their base morphisms $\underline \va$. 
Let us remark that $K$ is called  \emph{kernel} and geometries 
with trivial kernel are called \emph{effective}. If the geometry has non--trivial 
kernel, one can pass to the quotients $P/K \subset G/K$.

We are mainly interested in automorphisms of Cartan geometries.
It can be proved that the automorphism group $\Aut(\ba, \om)$  
of $(\ba \rightarrow M, \om)$ with $M$ connected is a Lie group of dimension 
at most dim$(G)$.
In particular, the automorphism group of connected components of 
the homogeneous model $(G\rightarrow G/P,\om_G)$ is exactly $G$.  
See \cite{S,parabook} for proofs. 
Considering these facts we suppose that $M$ is connected.

A \emph{parabolic geometry} is a Cartan geometry of 
type $(G,P)$ for a semisimple Lie group $G$ and  its parabolic 
subgroup $P$. The Lie algebra $\fg$ of the Lie group $G$ is then
equipped (up to the choice of Levi factor $\fg_0$ in $\fp$) with a grading 
of the form $\fg=\fg_{-k}\oplus \dots \oplus \fg_{0} \oplus \dots \oplus\fg_{k}$ 
such that the Lie algebra $\fp$ of $P$ is $\fp = \fg_0 \oplus \dots \oplus \fg_k$. 
We suppose that the grading of $\fg$ is fixed. 
By $G_0$ we denote the subgroup in $P$, with Lie algebra $\fg_0$,
consisting of all elements in $P$ whose adjoint action preserves the grading of $\fg$. 
Let us remark that $P$ is exactly the subgroup of elements of $G$ which preserve 
the usual filtration
$\fg=\fg^{-k}\supset \fg^{-k+1} \supset \dots\supset \fg^{k}=\fg_k$ of $\fg$ given by 
the grading, so $\fg^{i}=\fg_i \oplus \dots \oplus \fg_k$. 
A parabolic geometry corresponding to a grading of a length $k$ 
is called \emph{$|k|$--graded}. We are mainly interested in $|1|$--graded geometries 
and we formulate most of facts only for them.

The \emph{curvature} of a Cartan geometry is defined as 
$K:= d\om+\frac12[\om,\om]$, which is a strictly horizontal and 
$P$--equivariant two--form on $\ba$ with values in $\fg$.
It is fully described by a $P$--equivariant mapping $\ka: \ba \rightarrow 
\wedge^2 (\fg/\fp)^* \otimes \fg$, the so--called \emph{curvature function}.
If its values are in $\wedge^2 (\fg/\fp)^* \otimes \fp$, we call the geometry 
\emph{torsion--free}.
Notice that the Maurer--Cartan equation implies that the curvature of 
a homogeneous model vanishes. It can be proved that
if the curvature of a Cartan geometry vanishes, then it is
locally isomorphic to the homogeneous model of the same type, 
see \cite{S,parabook}.

For parabolic geometries, there is a notion on the normalization condition on the curvature. 
The $|1|$--graded geometry is called \emph{normal} if $\partial^* \circ \ka =0$, where 
$\partial^*:\wedge^2\fg_1\otimes\fg\rightarrow\fg_1\otimes\fg$ 
is given on decomposable elements by 
$\partial^*(X \wedge Y \otimes Z)=-Y\otimes[X,Z]+X\otimes[Y,Z]$
for $X,Y \in  \fg_{1} \simeq \fg_{-1}^*$ and $Z \in \fg$.

Let us remind that for each $|1|$--graded geometry $(\ba \rightarrow M,\om)$, there is a  
first order G--structure with a structure group~$G_0$ underlying this geometry:
The exponential mapping defines a diffeomorphism from $\fg_1$
onto a closed subgroup $P_+:=\exp \fg_1$ of $P$ and in fact, $P$ is a semidirect product 
of $G_0$ and the normal subgroup $P_+$, see \cite{parabook}. Because $P$ and $P_+$ act 
freely on $\ba$, we can form an orbit space $\ba_0:=\ba/ P_+$, which is a principal 
bundle $p_0: \ba_0 \rightarrow M$ with a structure group $P/P_+=G_0$. 
This is the reduction of a structure group of the frame bundle 
$\mathcal{P} := Gl(\fg_{-1}, TM)$  of $TM$ to the group $G_0$ with respect 
to the homomorphism $\Ad:G_0 \rightarrow Gl(\fg_{-1})$. (This is an isomorphism in the 
projective case, so $\ba_0$ is the full frame bundle in this case. 
The whole structure is given by the choice of a class of equivalent connections). 
Conversely, it can be proved the following fact, see \cite{parabook,AHS2} for details.
\begin{prop*}
Let $\ba_0 \rightarrow M$ be a reduction of the (first order) frame 
bundle $\mathcal{P}$ to the structure group $G_0$. 
Then there is a normal $|1|$--graded geometry $(\ba \rightarrow M,\om)$ 
of a suitable type inducing the given data.
\end{prop*}
Except for the projective structures, the $|1|$--graded normal geometry 
$(\ba \rightarrow M,\om)$ is unique up to an isomorphism and thus
there is the equivalence of the categories. 
In the projective case, there is an equivalence of categories 
between normal $|1|$--graded geometries 
of projective type and underlying frame bundles 
together with a class of projective equivalent torsion--free connections.

\subsection{Weyl structures and connections}		\label{1.2}
Let $(p:\ba\rightarrow M,\om)$ be a $|1|$--graded geometry of a type $(G,P)$ and 
let $\ba_0=\ba/P_+$ be the underlying bundle as above. We have 
the principal bundle $p_0:\ba_0 \rightarrow M$ with structure group $G_0$ and 
the principal bundle $\pi: \ba \rightarrow \ba_0$ with 
structure group $P_+$.
A \emph{Weyl structure} is a global smooth $G_0$--equivariant 
section $\si:\ba_0\rightarrow\ba$ of the projection~$\pi$.
For arbitrary $|1|$--graded geometries, Weyl structures always exist and any two 
Weyl structures $\si$ and $\hat \si$ differ by a $G_0$--equivariant mapping 
$\U:\ba_0\rightarrow \fg_1$ such that $\hat \si(u)= \si(u)\cdot \exp\U(u)$ 
for all $u \in \ba_0$. 
Since $\U$ gives the frame form of a 1--form on $M$, all Weyl 
structures form an affine space modeled on $\Om^1(M)$. 

In fact, any Weyl structure provides a reduction of the principal
bundle $\ba\rightarrow M$ to the subgroup $G_0\subset P$.
Given a Weyl structure $\si$, we can form the pullback
 $\si^*\om \in \Omega^1(\ba_0,\fg)$. This decomposes as 
$\si^*\om=\si^*\om_{-1}+\si^*\om_{0}+\si^*\om_{1}$
and $\si^*\om_{-1}$ is exactly the soldering form. 
The part $\si^*\om_0 \in \Omega^1(\ba_0,\fg_0)$ 
defines a principal connection on $p_0:\ba_0 \rightarrow M$ 
which we call a \emph{Weyl connection}.  
Moreover, this principal connection induces connections on all  
associated bundles.
In particular, we get a class of preferred affine connections
 on the tangent bundle characterized by the property 
that they share the same $\partial^*$--closed torsion.
We call each such  connection a Weyl connection, too, and we denote it by $\na$. 
The positive part $\si^*\om_1$ is called 
\emph{Rho--tensor} and is denoted by $\Rho$.

There are many formulas for the change of Weyl connections and 
related objects corresponding to the choice of various Weyl structures, which  
are analogous to the well known formulas from the conformal geometry. For instance,
the Rho--tensor transforms as 
$\hat \Rho(\xi) =\Rho(\xi)+\na_\xi \U +{1 \over 2}[\U,[\U,\xi]]$.
We will not need most of them explicitly, the reader can find them 
in \cite{parabook,TC,TBI}.

At the same time, the choice of $\si$ defines a decomposition of the curvature
$\si^*\ka = T+W+Y$ according to the values in $\fg_{-1}\oplus\fg_0\oplus\fg_1=\fg$.
The lowest part $T$ of the decomposition corresponds to the torsion of $\na$. 
It does not depend on the choice of the Weyl structure and thus it is an invariant of 
the parabolic geometry. In fact, it coincides with the Cartan torsion.
The part $W$ is called  \emph{Weyl curvature}. It can be written via the Lie 
algebra differential $\partial$ as $W=R+\partial \Rho$, 
where $R$ is the curvature of $\na$, see \cite{parabook, TBI} for details.
If the torsion of the geometry vanishes, then the Weyl curvature 
is independent of the choice of the Weyl structure and is an invariant 
of the parabolic geometry. The positive part $Y$ is called 
\emph{Cotton--York tensor}.

\subsection{Normal Weyl structures} \label{1.3}
Among general Weyl structures, there are various specific subclasses. 
From our point of view, the most interesting are so called 
normal Weyl structures. 
Denote flows of constant vector fields 
$\omi(X) \in \frak X(\ba)$ as $\Fl^{\omi(X)}_t(u)$.
We can define a \emph{normal Weyl structure} at $u$ 
as the only $G_0$--equivariant section $\si_u:\ba_0 \rightarrow  \ba$ 
with the property
$
\si_u \circ \pi \circ {\rm Fl}_1^{\omi(X)}(u) = {\rm Fl}_1^{\omi(X)}(u)
.$ 
Although the normal Weyl structure is indexed by $u\in\ba$, 
it clearly depends only on the $G_0$--orbit of $u$.
In general, normal Weyl structures are defined locally over some neighborhood 
of $p(u)=:x$. 
They are closely related 
to the normal coordinate systems for parabolic geometries and generalize the affine 
normal coordinate systems, see \cite{parabook,CSZ}.
There is also another useful characterization of normal Weyl structures: 
The Rho--tensor $\Rho$ of the normal Weyl structure at $u$ has the property
that for corresponding Weyl connection $\na$ and for all $k\in\Bbb N$, 
the symmetrization of $(\na_{\xi_k}\dots\na_{\xi_1}\Rho)(\xi_0)$ over all $\xi_i\in TM$ 
vanishes at $x$. In particular, $\Rho(x)=0$, see \cite{parabook,WS} for proofs.

\section{Basic facts on symmetries} \label{2}
In this section, we remind the definition and summarize quickly properties of
symmetries on $|1|$--graded parabolic geometries.  
For a complete discussion and proofs see \cite{ja-srni06,ja-elsevier} and 
references therein. 

\subsection{Definitions} \label{2.1} 
Let $(\ba \rightarrow M, \om)$ be a $|1|$--graded parabolic geometry 
of a type $(G,P)$. A \emph{symmetry} centered at $x$ 
is a diffeomorphism $s_x$ of $M$ such that:
\begin{itemize} 
\item[(i)] $s_x(x)=x$, 
\item[(ii)] $T_x s_x =-\operatorname{id}$ on $T_xM$, 
\item[(iii)] $s_x$ is covered by an automorphism $\va_x$ of the parabolic geometry. 
\end{itemize}
We call the parabolic geometry \emph{symmetric} if there is a symmetry at 
each point $x \in M$. 
It can be proved that each $s_x$ is involutive and $x$ is its isolated fixed point, see
\cite{ja-elsevier}.
Under the assumption that the geometry is effective, 
the (uniquely given) covering $\va_x$ of $s_x$ is also involutive. 

From the usual point of view, the definition above reflects the classical 
notion of affine symmetries. The first two properties simply say that 
our symmetries follow the classical intuitive idea. The third one means that
the latter symmetries respect the underlying geometrical structure given by the
$|1|$--graded geometry. 

\subsection{Flat models} \label{2.2-3}
The simplest candidates for symmetric geometries are homogeneous models 
$(G \rightarrow G/P, \om_G)$.
It is well known that all automorphisms of (connected components of) 
homogeneous models are exactly left multiplications by 
elements of $G$, see \cite{S, parabook}.
Transitivity of the action of $G$ on $G/P$ says that to decide whether the 
homogeneous model is symmetric, it suffices to find an element of $G$ giving 
a symmetry at the origin.
We have the following statement, see \cite{ja-srni06}:
\begin{prop*}
All symmetries of the homogeneous model centered at the origin are parametrized 
by elements $g_0 \exp Z \in P$, where $Z\in\fg_{1}$ is arbitrary and 
$g_0 \in G_0$ is such that $\Ad_{g_0}=-\rm{id}$ on $\fg_{-1}$. 
\end{prop*} 

Because we work with effective geometries, there is at most one element 
$g_0 \in G_0$ with the  above property  and it is usually a simple exercise to find it.
If the element exists, then the homogeneous model is symmetric and there is an 
infinite amount of symmetries at each point. If there is no 
such element, then no Cartan geometry of the same type is symmetric. 
As an example, we show here the projective case. Analogous computations showing 
that the homogeneous models of many types of $|1|$--graded geometries are symmetric
can be found in \cite{ja-dis}.

\begin{exam}\emph{A homogeneous model of a non--oriented projective geometry.} \label{exam-proj}
Consider $G=PGl(m+1,\mathbb{R})$, the quotient of
$Gl(m+1,\mathbb{R})$ by the subgroup of all multiples of the identity. 
This group acts on the space of lines passing through the origin in $\R^{m+1}$ 
and as $P$ we take the stabilizer of the line generated by the first standard basis 
vector. Clearly $G/P \simeq \mathbb{R} P^m$. 
With this choice, $(G \rightarrow G/P,\om_G)$ is the homogeneous model
of an effective $|1|$--graded geometry. The elements of the Lie algebra
$\fg = \frak{sl}(m+1,\mathbb{R})$ are of the block form $\left( \begin{smallmatrix}
-tr(A)&Z\\X& A\end{smallmatrix} \right)$ with blocks of sizes $1$ and $m$, 
where $X \in \mathbb{R}^m \simeq \fg_{-1}$, $Z \in \mathbb{R}^{m*} \simeq \fg_{1}$
and $\fg_0$ is the block--diagonal part which is determined by 
$A \in \frak{gl}(m,\mathbb{R})$.
The subgroup $G_0$ is isomorphic to $Gl(m,\mathbb{R})$ because each class
in $G_0$ has exactly one representative of the form 
$\left( \begin{smallmatrix} 1&0\\0& B \end{smallmatrix} \right)$.

For each $a \in G_0$ represented by some $\left( \begin{smallmatrix}
1&0\\0& B\end{smallmatrix} \right)$ 
and for $V=\left( \begin{smallmatrix}
0&0\\X& 0\end{smallmatrix} \right) \in \fg_{-1}$, 
the adjoint action $\operatorname{Ad}_aV$ is  given by $X \mapsto BX$ and
we look for elements $\left( \begin{smallmatrix}
1&0\\0& B\end{smallmatrix} \right)$ such that $BX=-X$ for each 
$X$. It is easy to see that we may represent the only prospective solution 
by 
$\left( \begin{smallmatrix}
1&0\\0& -E\end{smallmatrix} \right)$.
This element represents a class in $G_0$
and yields a symmetry. All elements giving some symmetry 
at the origin are then represented by matrices of the form 
$\left( \begin{smallmatrix} 1&W\\0& -E\end{smallmatrix} \right)$ 
for all $W \in \mathbb{R}^{m*}$. Thus homogeneous models of 
non--oriented projective geometries are symmetric. 
\end{exam}

\subsection{Local properties of symmetries} \label{2.2}
Let us mention here that many properties of symmetries 
on $|1|$--graded geometries are similar to the classical symmetries. 
First recall that an invariant tensor field of odd degree which is 
invariant with respect to a symmetry at $x\in M$ vanishes at $x$, see \cite{ja-srni06}. 
This applies directly to the Cartan torsion, which 
is an invariant tensor of degree three.
Thus symmetric $|1|$--graded geometries are torsion--free.
Then, in particular, Weyl connections are torsion--free.

Let $s_x$ be a symmetry at $x$. Denote $\va$ its covering and $\va_0$ the 
corresponding underlying morphism. 
An arbitrary Weyl structure $\si$ satisfies 
$ 
\va(\si(u_0))= \si(\va_0(u_0))\cdot \exp F(u_0)
$
for a suitable $G_0$--equivariant function $F:\ba_0 \rightarrow \fg_1$. 
In fact, there exists a Weyl structure $\hat \si$ such that 
\begin{align} \label{jedna} 
\va(\hat \si(u_0))=\hat \si(\va_0(u_0))
\end{align} 
holds in the fiber over $x$.
In the fiber over $x$, $\hat \si$ is given uniquely (although there can exist 
various Weyl structures with this property).
Moreover, there is at least one Weyl structure satisfying $(\ref{jedna})$
over some neighborhood of $x$. This property has the normal Weyl structure 
$\si_u$ given at $u=\hat \si(u_0)$, $u_0 \in p_0^{-1}(x)$  
for an arbitrary $\hat \si$ satisfying $(\ref{jedna})$ in the fiber over $x$.
See \cite{ja-dis, ja-elsevier} for proofs. 
We will also return to this fact in more detail in \ref{4.2}.

These observations allow us to use the following terminology: 
For each symmetry $s_x$ at $x$ on an effective $|1|$--graded geometry 
and for its (uniquely given) covering $\va$, we call the latter Weyl structures 
\emph{invariant with respect to $s_x$} or shortly \emph{$s_x$--invariant}
at the point $x$ or on a neighborhood of $x$, respectively. 
In these terms, the above facts can be formulated as follows (see \cite{ja-elsevier}):  
\begin{prop*}
Let $s_x$ be a symmetry at $x$ on a $|1|$--graded geometry. 
There exist Weyl structures which are invariant with respect to $s_x$ 
in the fiber over $x$ and all of them coincide in the fiber over $x$. 
Moreover, at least one of them is invariant with respect to $s_x$ over some 
neighborhood of $x$.
\end{prop*}
Suppose there is a symmetry $s_x$ at $x$ on a $|1|$--graded geometry. 
Consider a Weyl connection $\na$ corresponding to a Weyl structure 
which is invariant on some 
neighborhood of $x$ and denote $T,R,W$, and $\Rho$ its torsion,
curvature, Weyl curvature and Rho--tensor, respectively.
Because the Cartan torsion vanishes at $x$, we have $T=0$ at $x$. 
Then $W$ is an invariant tensor (of degree four) at $x$ and consequently,
$\na W$ is an invariant tensor of degree five and thus 
$\na W=0$ at $x$. Because $\na$ is invariant locally around $x$ 
(and not only at the point $x$), $\Rho$ is invariant at $x$.
The formula $W = R +\partial \Rho$ then gives that $R$ has to be an invariant 
of degree four and consequently $\na R=0$ at $x$.

\begin{cor*}
Suppose there is a symmetry $s_x$ centered at $x$ on a $|1|$--graded 
geometry. On a neighborhood of $x$, there exists a Weyl
connection $\na$ invariant under $s_x$. In particular, $T=0$ and $\na R=0$ at $x$. 
\end{cor*}
In the other words, each symmetry is (locally) an affine symmetry for some admissible 
affine connection, e.g. the Weyl connection corresponding to the invariant normal 
Weyl structure. For this normal connection, the symmetry is easily 
understandable: It is simply reverting of geodesics of the invariant normal 
connection (which are also generalized geodesics of the Cartan connection). 
See \cite{CSZ} and \cite{ja-elsevier} for details.
Contrary to the affine case, we know nothing about an invariance of possible 
invariant Weyl structures with respect to various 
symmetries at different points.

\subsection{Non--flat symmetric geometries} \label{2.3}
As we have seen in \ref{2.2-3}, there can exist many different symmetries at one point. 
But the existence of a non--zero curvature gives quite strong restrictions on the 
number of different symmetries at one point.
Let us summarize here shortly restrictions given on the number of possible symmetries 
for non--flat geometries. See \cite{ja-srni06, ja-elsevier} for all details.

Let us remind firstly, that it is not difficult to find all $|1|$--graded geometries 
(with simple Lie group $G$). It corresponds to the well known classification 
of semisimple Lie algebras and their parabolic subalgebras, see \cite{Y, parabook}.
See also \cite{WS,ja-dis} for the list of them.
In fact, using the classification it can be proved that most types of normal
$|1|$--graded geometries have to be locally isomorphic to the homogeneous model, 
if they are symmetric, see \cite{ja-srni06}. (Here the normality is only some technical 
restriction which plays no role, if we understand symmetries as morphisms of 
the underlying G--structure.)

Only a few types of $|1|$--graded geometries can carry a symmetry at the point 
with a non--zero curvature. 
We summarize the list of them with corresponding restrictions in the following 
statement, see \cite{ja-elsevier, ja-dga07} for proofs.
\begin{prop*}
(1) There can exist at most one symmetry at the point with a non--zero curvature on 
projective geometries, almost quaternionic geometries and conformal geometries 
of positive definite or negative definite signature.
\\
(2) Suppose there are two different symmetries with the center at $x$ 
on an indefinite conformal geometry or on the almost Grassmannian geometry of type 
$(2,m)$ or $(m,2)$, respectively, and denote $\si_1, \si_2$ their invariant Weyl 
structures. Suppose that the function $\U$ given by $\si_1=\si_2 \cdot \exp \U$
has a non--zero length or a maximal rank at $x$, respectively. 
Then the Cartan curvature vanishes at $x$.
\end{prop*}
It is not surprising that these are exactly the best known types of $|1|$--graded 
parabolic geometries. 
From now, we discuss mainly these examples because they are also the most 
interesting ones from our point of view. 

Let us return briefly to the condition in the part (2) of the Proposition.
First, the function $\U$ is correctly defined at the point $x$, because
all invariant Weyl structures of an arbitrary fixed symmetry coincide at $x$.
Moreover, it is not difficult to see that different symmetries have different 
invariant Weyl structures. Clearly, each symmetry determines uniquely the invariant 
normal Weyl structure. If these Weyl structures were equal for two symmetries, 
then the two symmetries would be equal as geodesical symmetries, i.e. they 
would give the same reverting of geodesics of the invariant normal Weyl connection.
One can also see that the result does not 
depend on the order of the two symmetries. Thus we get that there could exist 
more symmetries at one point, but the non--zero difference $\U$ of each two of them 
is degenerate in some sense.

\section{Systems of symmetries} \label{3}
In this section, we start the discussion of systems of symmetries on 
$|1|$--graded geometries. In the view of the known fact that each affine 
symmetric space is homogeneous and its system of symmetries is smooth, 
we concentrate first on smooth systems of symmetries on an arbitrary 
$|1|$--graded geometry.

\subsection{Smooth systems of symmetries} \label{3.2}
One of the usual definitions says that a homogeneous space is a manifold $M$ 
with transitive action of a Lie group $H$. Suppose there is a parabolic
geometry $(\ba \rightarrow M,\om)$ defined on this $M$. 
We call it \emph{homogeneous parabolic geometry}, if $H$ acts transitively 
by automorphisms of the parabolic geometry. More precisely, we require that 
each element of $H$ (viewed as an automorphism of $M$) has covering in 
$\Aut(\ba,\om)$. In fact, it is irrelevant whether we understand the Lie 
group $H$ as a subgroup of automorphisms of $M$ or subgroup of their coverings. 
We will not distinguish between them and we will simply write  
$H \subset \Aut(\ba,\om)$. 

Let us now discuss smooth systems of symmetries on an arbitrary $|1|$--graded geometry 
$(\ba \rightarrow M,\om)$, i.e. suppose we have a symmetry $s_x$ fixed at each 
$x$ such that the map $S:M \times M \rightarrow M$ given by $S(x,y)=s_x(y)$ is smooth.
\begin{prop*}
Let $(\ba \rightarrow M, \om)$ be a symmetric $|1|$--graded geometry 
with $M$ connected and suppose that the system of symmetries 
$S: M \times M \rightarrow M,\ (x,y) \mapsto s_x(y)$ is smooth. 
Then the geometry is homogeneous. More precisely, a subgroup 
$H \subset \Aut(\ba, \om)$ containing all symmetries 
acts transitively on $M$.
\end{prop*}
\begin{proof}
Choose an arbitrary $x_0 \in M$ and discuss the map $f:M \rightarrow M,\ x \mapsto s_x(x_0)$. 
The map $f$ is smooth and preserves $x_0$ because $f(x_0)=s_{x_0}(x_0)=x_0$. 
Thus we can compute $T_{x_0}f:T_{x_0}M \rightarrow T_{x_0}M$. Take a curve $c:I\rightarrow M,\ 0 \mapsto x_0$ representing some $\xi(x_0)\in T_{x_0}M$ as $\xi(x_0)= {d \over dt}|_{0} c=c'(0)$. 
Then the curve $I \rightarrow M$, $ t \mapsto s_{c(t)}(x_0)$ represents $Tf.\xi(x_0)$ as 
${d \over dt}|_{0}s_{c(t)}(x_0)$. 
We compute ${d \over dt}|_{0}s_{c(t)}(x_0)$ using the obvious identity $s_{c(t)}(c(t))=c(t)$. 
Differentiating at $t=0$ gives that
$$
{d \over dt}|_{0}s_{c(t)}(c(t))= {d \over dt}|_{0}s_{c(t)}(x_0)+ T_{x_0}s_{x_0}.c'(0)=
{d \over dt}|_{0}s_{c(t)}(x_0)-c'(0)
$$
equals to $c'(0)$ and we get
$$
{d \over dt}|_{0}s_{c(t)}(x_0)=2 \cdot c'(0).
$$

Thus $T_{x_0}f=2\cdot \rm{id}$ and $f$ is locally invertible in some neighborhood of $x_0$.
This shows that the orbit $N$ of $x_0$ is open (with respect to symmetries 
and thus under the action of the automorphism group). 
The complement of $N$ is also open because 
it is a union of orbits of points from $M \setminus N$ and these orbits are open 
from the same reason. Thus $N$ is simultaneously closed and we get $M=N$ from the 
connectedness of $M$. All together, $H \subset \Aut(\ba,\om)$ containing symmetries 
acts transitively and thus $(\ba\rightarrow M,\om)$ is a homogeneous $|1|$--graded geometry.
\end{proof}

Thus in the case of a smooth system of symmetries on $(\ba \rightarrow M, \om)$, 
we can write $M=H/K$ for suitable Lie groups $H$ and $K$. More precisely, $H$ is as above, 
$K \subset H$ is a stabilizer of the point $x_0 \in M$ and the isomorphism 
is given by the mapping $H \rightarrow M$, $h \mapsto h(x_0)$ which factorizes correctly
to $H/K$  because $H$ acts transitively on $M$. 
The Proposition also says that if a symmetric $|1|$--graded geometry is not 
homogeneous, then there is no smooth system of symmetries given on the geometry.

\subsection{A non--trivial example} \label{3.3}
Contrary to the affine case, there exist symmetric spaces which are not 
homogeneous and thus the above observations are not trivial. 
An idea of a construction of such projectively symmetric examples was given in \cite{P} 
and we sketch the example in the language of parabolic geometries.

\begin{exam} \emph{A non--homogeneous projectively symmetric space.} \label{exam-nehomog}
We start with the homogeneous model for projective structures, 
see Example \ref{exam-proj}. 
Thus take $(G \rightarrow G/P,\om_G)$ where $G=PGl(m+1,\R)$ and if we denote 
$e_1,\dots,e_{m+1}$ the standard basis, then $P$ is the stabilizer of 
the line generated by $e_1$. 
Clearly, the automorphism group $G$ acts transitively on $G/P\simeq\R P^m$. 
Let us define a new manifold $M$ such that we remove two points from the model 
$G/P \simeq \R P^{m}$, e.g. points given by the last two basis vectors $e_{m}$ and $e_{m+1}$. 
One can easily reduce the projective structure from $G/P$ to $M$ and we get a 
$|1|$--graded geometry on $M$.

First, we show that the automorphism group 
does not act transitively on $M$. One can easily describe this group: 
Its elements correspond to those automorphisms from $G$ which maps $M$ to $M$, 
i.e. which either preserve  the points of 
$G/P$ given by $e_m$ and $e_{m+1}$ or which swap these two points.
We can represent them by matrices from $G$ with last two columns of the form
$\left( \begin{smallmatrix} 0 &\dots &0 & *& 0 
\\ 0&\dots & 0 & 0 &*\end{smallmatrix} \right)^T$ 
or
$\left( \begin{smallmatrix} 0 &\dots &0 & 0& *  
\\ 0&\dots & 0 & * &0 \end{smallmatrix} \right)^T$. 
Exactly these automorphisms of the homogeneous model can be correctly restricted to 
automorphisms of $M$. 
Choose an arbitrary point of $M$ which lies on the projective line (on $G/P$)
generated by points given by $e_m$ and $e_{m+1}$. It is represented by a vector 
$w=\left( \begin{smallmatrix} 
0& \cdots & 0&x_m&x_{m+1}
\end{smallmatrix} \right)^T$ 
where $x_m$ and $x_{m+1}$ are both nonzero. (If one of them is zero, 
then we have exactly some removed point.) 
One can see from the matrices that there is no automorphism 
mapping $w$ out of the line. 

Second, let us show that $M$ is symmetric. There are many symmetries 
at the origin $o \in M \subset G/P$. We know that all symmetries of the 
homogeneous model at $o$ are given by elements of the form 
$s_o= \left( \begin{smallmatrix} 1&W\\0& -E\end{smallmatrix} \right)$, see Example \ref{2.2-3}.
To get a correctly defined symmetry on $M$, it suffices to take an arbitrary element  
with last two coordinates of $W$ vanishing. Denote an arbitrary such element by $s$. 
Clearly, $s$ preserves the removed points. 
Now it is easy to get a symmetry at each $x$ which does not lie on the  
projective line described before.
Such $x$ is represented by some 
$v=\left( \begin{smallmatrix} 
x_1 & \dots & x_{m-1} & x_{m} & x_{m+1}
\end{smallmatrix} \right)^T$ 
where at least one of $x_1, \dots, x_{m-1}$ is non--zero.
Then elementary linear algebra says that we can find $g \in G $, $g \cdot e_1=v$ 
such that $g$ preserves the removed points. We take a symmetry at $x=gP$ of the form 
$s_x=g  s g^{-1}$.

Points $x$ on the line, represented by some vectors
$w=\left( \begin{smallmatrix} 0&\dots&0&x_m&x_{m+1}\end{smallmatrix} \right)$ as above,  
are more complicated. There are $g \in G$ such that 
$g \cdot e_1=w$, but non of them preserves or swaps $e_m$ and $e_{m+1}$. Thus take
an arbitrary fixed such $g$, compute  expression $g s_og^{-1}$ 
(here $s_o=\left( \begin{smallmatrix} 1&W\\0& -E\end{smallmatrix} \right)$ as above) 
and analyze the result depending on $W$.
For example, choose $g$ such that $e_1 \mapsto w$, $e_2\mapsto e_2,$ $\dots$, $e_m\mapsto e_m$,
$e_{m+1}\mapsto e_1$.
If we compute $g s_og^{-1}$ (in the standard basis), its last two columns are of the form
$\left( \begin{smallmatrix} 
0 &\dots &0 & x_m v_m -1& x_{m+1}v_m \\ 
0&\dots & 0 &{x_m \over x_{m+1}}(x_mv_m+2) & 1-x_m v_m 
\end{smallmatrix} \right)^T$, 
where we denote $W=\left( \begin{smallmatrix} v_2& \dots& v_{m+1}\end{smallmatrix} \right)$. 
We see from the matrix that for $v_m={1 \over x_m}$, the composition 
is a symmetry which swaps the removed points and no symmetry preserves them.
\end{exam}

Consequently, there is no smooth system of symmetries on the geometry, although there
can exist many symmetries at each point. 
Let us remark, that the latter space is locally isomorphic to the homogeneous model, 
so there is no contradiction to Proposition \ref{2.3}.

\subsection{Homogeneous parabolic geometries} \label{3.1}
First note that if $(\ba\rightarrow M,\om)$ is homogeneous, than $M=H/K$ for
suitable $H\subset \Aut(\ba,\om)$ acting transitively on $M$ and $K\subset H$ is a
stabilizer of a point $x_0 \in M$. There exists a structure
of a smooth manifold on $M$ such that the canonical projection $q:H \rightarrow H/K$ 
is a surjective submersion, see \cite{KMS}. 
We can prove a partial converses to the Proposition \ref{3.2}. 

\begin{prop*}
Let $(\ba \rightarrow M, \om)$ be a $|1|$--graded parabolic geometry and suppose that 
$H \subset \Aut(\ba,\om)$ acts transitively on $M$. Suppose further 
that for some $x_0 \in M$ there is a symmetry $s_{x_0}$ such that 
$k \circ s_{x_0} \circ k^{-1}=s_{x_0}$  for all $k \in H$ such that 
$k(x_0)=x_0$. Then there exists a smooth system of symmetries $S$.
\end{prop*}
\begin{proof}
Let $s_{x_o}$ be the symmetry at $x_0 \in M$.
Define a map $f: H \rightarrow \Aut(\ba, \om)$ by 
$f(h)=h \circ s_{x_0} \circ h^{-1}$.
The map $f$ is correctly defined and smooth because it is simply
multiplication of three elements from a Lie group $\Aut(\ba, \omega)$.

We show that $f$ factorizes correctly to the mapping 
$g:H/K \rightarrow \Aut(\ba, \om)$
via the surjective submersion $q:H \rightarrow H/K$,
where $K \subset H$ consists of elements which preserve $x_0 \in M$.
For $h \in H$ we have 
$$
g(hK)= f(h)= h \circ s_{x_0} \circ h^{-1}
.$$
If we change a representative of the class $hK$, i.e. we take 
$h \circ k$ for some $k \in K$, we get
$$
g(hkK)=f(h \circ k)= h \circ k \circ s_{x_0} \circ k^{-1} \circ h^{-1}.
$$
The assumption  $k \circ s_{x_0}\circ k^{-1}=s_{x_0}$ gives 
$$
h \circ k \circ s_{x_0} \circ k^{-1} \circ h^{-1}=h \circ s_{x_0} \circ h^{-1}
$$ 
and the factorization $f=g \circ q$ is correctly defined.
The universal property of the surjective submersion $q$ says that the map 
$g:H/K \rightarrow \Aut(\ba, \om)$ is smooth. 

Transitivity of the action $H$ on $M$ says that $M \simeq H/K$. 
Thus we have defined correctly a symmetry $s_x$ at each 
$x \in M$ such that $s_x = h \circ s_{x_0} \circ h^{-1}$ for $x=h(x_0)$
and such that the mapping $x \mapsto s_x$ is smooth. Equivalently, 
$S:M \times M \rightarrow M$, $(x,y) \mapsto s_x(y)$ is a smooth system of symmetries.
\end{proof}
In particular, if there is exactly one symmetry $s_x$ at each $x$, then 
$s_x = h \circ s_{x_0} \circ h^{-1}$ is clearly satisfied for each $x=h(x_0)$ and
we have the following fact:
\begin{cor*}
Suppose there exists exactly one symmetry at each point on a homogeneous 
$|1|$--graded geometry. Then the system of symmetries is smooth.
\end{cor*}

\section{The invariant connection} 
In this section, we continue in the discussion of systems of symmetries on arbitrary 
$|1|$--graded geometries. Motivated by \cite{P}, we show here that under quite 
weak conditions on the system, there exist an admissible affine connection, 
which is invariant with respect to all symmetries from the system. 

\subsection{Fiberwise invariant Weyl structures} \label{4.1}
In Section \ref{2}, we discussed each symmetry on a $|1|$--graded geometry 
separately. We know that for each $s_x$, there exist Weyl 
structures which are invariant with respect to $s_x$ in the fiber over $x$ and 
that some of them are invariant on some neighborhood of $x$, 
see also \cite{ja-elsevier}. 
Suppose there is a system of symmetries $S: M \times M\rightarrow M$, $S(x,y)=s_x(y)$ 
defined on a $|1|$--graded geometry. 
There are natural questions whether some of the Weyl structures are invariant with 
respect to all symmetries from the system $S$ and how many such Weyl structures can exist?
We start the discussion with the following observation:

\begin{lemma*}
Let $S: M \times M\rightarrow M$, $S(x,y)=s_x(y)$ be an arbitrary system of 
symmetries on a $|1|$--graded geometry $(\ba \rightarrow M,\om)$. There can exist 
at most one Weyl structure 
$\hat \si$ which is invariant with respect to all symmetries from $S$.
\end{lemma*}
\begin{proof}
If there exists a Weyl structure $\hat \si$ which is invariant with respect to all 
symmetries from the system $S$, then, in particular, $\hat \si$ is invariant with 
respect to $s_x$ in the fiber over $x$. This holds for each $x \in M$ and 
we use this fact to construct the only possible candidate for an invariant Weyl structure.

Denote by $\si_x$ an arbitrary fixed $s_x$--invariant Weyl structure at $x$.
Define a mapping 
\begin{align*}
\hat \si :\ &\ba_0 \rightarrow \ba  
\\& u_0 \mapsto \si_x(u_0)\ \ \ {\rm for} \ \ \ p_0(u_0)=x. 
\end{align*}
In other words, we glue together the invariant Weyl structures in appropriate fibers. 
Because all $s_x$--invariant Weyl structures coincide in the fiber over $x$, the definition 
of the mapping $\hat \si$ is correct.
All $\si_x$ are $G_0$--equivariant and thus $\hat \si$ is a $G_0$--equivariant mapping 
$\ba_0\rightarrow \ba$. If this mapping is smooth, 
then it is a Weyl structure. It is clear from the construction that $\hat \si$ is 
the only possible candidate for an invariant Weyl structure.
\end{proof}

It is not necessarily true in general that the mapping $\hat \si$ constructed 
in the Lemma for some system $S$ is smooth. For instance, the geometry from 
Example \ref{exam-nehomog} cannot carry a system of symmetries with corresponding 
$\hat \si$ smooth although there exist many symmetries at each point and thus for any 
system of symmetries, there is no invariant Weyl structure. Conversely, $\hat \si$ is 
smooth and thus a Weyl structure for each smooth system $S$. 
\begin{def*}
If the mapping $\hat \si: \ba_0 \rightarrow \ba$ induced by the system $S$ as above 
is smooth, we call the section $\hat \si$ the \emph{fiberwise invariant Weyl structure} 
of the system~$S$. In this case, we say that the system $S$ admits fiberwise invariant Weyl 
structure.
\end{def*}

The motivation for the name is clear. Such a Weyl structure $\hat \si$ still need not to 
be invariant. But for all $x \in M$, it is invariant with respect to $s_x$ in the 
fiber over $x$ (which is a much weaker property). However, it is still quite interesting 
object. For instance, there is the following statement:
\begin{prop*}
Suppose that the system $S$ of symmetries admits fiberwise invariant Weyl structure 
on a $|1|$--graded geometry. Then the corresponding Weyl connection $\na$ satisfies
$T=0$ and $\na W=0$.
\end{prop*}
\begin{proof}
Because our geometry is symmetric, the Cartan torsion vanishes and thus 
each Weyl connection is torsion--free.
Then the Weyl curvature $W$ is an invariant of the geometry. 
For each $x$, the fiberwise invariant Weyl 
structure $\hat \si$ is $s_x$--invariant in the fiber over $x$. Then for the
corresponding Weyl connection we have that $\na W$ is an invariant tensor of 
degree five at each $x \in M$ and thus $\na W=0$.
\end{proof}
Let us remark that this does not imply, that the curvature $R$ of $\na$ is invariant and 
covariantly constant. The difference between invariance at the point and invariance 
on the neighborhood is crucial here. One can see from the formula for the change of 
Rho--tensor (see \ref{2.2}) that it depends on the derivative of the change. Invariance 
at the point is not sufficient and fiberwise invariant Weyl structure is too weak.

\subsection{Invariant Weyl structures} \label{4.2}
For a suitable system of symmetries, the fiberwise invariant Weyl structure is the 
only Weyl structure which could be invariant with respect to all symmetries. 
We would like to know when this Weyl structure really is invariant.
\begin{prop*}
Let $(\ba \rightarrow M,\om)$ be a $|1|$--graded geometry and 
let $S: M \times M\rightarrow M$, $S(x,y)=s_x(y)$ be a system of symmetries. 
There exists a Weyl structure which is invariant with respect to all symmetries from $S$ 
if and only if 
\begin{enumerate}
\item[(i)] the system $S$ admits fiberwise invariant Weyl structure,
\item[(ii)] $s_x \circ s_y \circ s_x = s_{s_x(y)}$ holds for each $x,y \in M$. 
\end{enumerate}
\end{prop*}
\begin{proof}
Let us start with some notation.
For the system of symmetries $S:M \times M \rightarrow M$, $S(x,y)=s_x(y)$, 
there is the usual covering of the mapping $S$ of the form 
$\va: M \times \ba \rightarrow \ba$, i.e.
$\va(x, \cdot)$ is simply the covering of $s_x$. 
Clearly, there is also the corresponding underlying morphism 
$\va_0: M \times \ba_0 \rightarrow \ba_0$.
For each $x \in M$, $u_0 \in \ba_0$ and an arbitrary Weyl structure $\si$, we can write
\begin{align} \label{dva}
\va(x,\si(u_0))= \si(\va_0(x,u_0))\cdot \exp F(x,u_0)
\end{align}
for a suitable function $F:M \times \ba_0 \rightarrow \fg_1$.
In this notation, the Weyl structure $\si$ is invariant if and only if 
$$ 
\va(x,\si(u_0))=\si(\va_0(x,u_0))
$$
holds for each $x \in M$ and $u_0 \in \ba_0$. 
Let us remark, that the fiberwise invariant Weyl structure in general 
satisfies this equality only 
at the points $(p_0(u_0),u_0) \in M \times \ba_0$. 
\smallskip

Suppose that the Weyl structure $\hat \si$ is invariant with respect to $S$. 
Then $\hat \si$ is clearly fiberwise invariant Weyl structure for the system 
$S$ and we get (i). 
Now, for (coverings of) $s_x$ and $s_y \circ s_x$ 
we have
\begin{align*}
\va(x,\hat \si(u_0))&=\hat \si (\va_0(x,u_0)),\\
\va(y,\va(x,\hat \si(u_0))&=\va(y,\hat \si(\va_0(x,u_0)))= \hat \si(\va_0(y,\va_0(x,u_0)))
\end{align*} thanks to the invariance of $\hat \si$ with respect to $s_x, s_y \in S$.
Consequently, for (the covering of) $s_x \circ s_y \circ s_x$ we get 
$$
\hat \si(x,\va(y,\va(x,\hat \si(u_0))))=\va(x,\hat \si(\va_0(y,\va_0(x,u_0))))=
\hat \si(\va_0(x,\va_0(y,\va_0(p,u_0))))
$$
for each $u_0 \in \ba_0$ and $x,y \in M$. In other words, 
the Weyl structure $\hat \si$ is invariant with
respect to the composition $s_x \circ s_y \circ s_x$.
Let us shortly verify here, that the mapping $s_x \circ s_y \circ s_x$ satisfies  
conditions of a symmetry: 
Clearly, $(s_x \circ s_y \circ s_x)(s_x(y))=s_x(y)$ from the involutivity of $s_x$
and the definition of $s_y$. Thus $s_x(y)$ is its fixed point.
Then for $\xi(s_x(y)) \in T_{s_x(y)}M$ we have that
$T(s_x \circ s_y \circ s_x).\xi(s_x(y))$
is equal to
$$
T_y s_x.T_y s_y.T_{s_x(y)}s_x.\xi(s_x(y))=
T_y s_x.(-T_{s_x(y)}s_x.\xi(s_x(y)))=-\xi(s_x(y))
.$$ 
Thus the composition is some symmetry at the point $s_x(y)$ and we 
know that the Weyl structure $\hat \si$ is invariant with respect to this symmetry. 
But $\hat \si$ is also invariant with respect to the symmetry $s_{s_x(y)}$ 
from the system $S$. 
If two symmetries share the same invariant Weyl structure, they are equal, see  \ref{2.3}, 
and we get $s_x \circ s_y \circ s_x = s_{s_x(y)}$ for each $x,y \in M$ which is (ii).
\smallskip

Conversely, we show that conditions (i) and (ii) imply the existence of an invariant 
Weyl structure or, equivalently, that under the condition (ii), the fiberwise 
invariant Weyl structure given by $S$ is really invariant with respect to $S$.

Let $\si$ be an arbitrary fixed Weyl structure. 
First, we find a useful condition on $\U$ such that 
$\hat\si(u_0)=\si(u_0) \cdot \exp \U(u_0)$
is the invariant Weyl structure. 
Substitute the formula for the change of Weyl structures into the 
expression $(\ref{dva})$. The left hand side is then equal to
\begin{align*}
\va(x,\hat \si(u_0)) &=  \va(x,\si(u_0)) \cdot \exp \U(u_0) \\
& =  \si(\va_0(x,u_0))\cdot \exp F(x,u_0) \cdot \exp \U(u_0) 
\end{align*}
and the right hand side is of the form
\begin{align*}
\hat\si(\va_0(x,u_0))=\si(\va_0(x,u_0))\cdot \exp \U(\va_0(x,u_0)).
\end{align*}
Comparing of these expressions gives the equation
\begin{align} \label{tri}
F(x,u_0)= \U(\va_0(x,u_0)) - \U(u_0),
\end{align}
where the unknown quantity is the change $\U$, while $F$ is given by $\si$. 
Clearly, $\hat \si=\si \cdot \exp \U$ is invariant if and only if $\U$ solves 
$(\ref{tri})$ on all of $M \times \ba_0$. Let us remark, that $\hat \si$ is 
fiberwise invariant if and only if $\U$ satisfies the equation only on the points 
of the form $(p_0(u_0),u_0) \in M \times \ba_0$. 

Let us discuss the equation $(\ref{tri})$ in more detail.
It is not difficult to find $\U$ such that $\hat \si=\si \cdot \exp \U$ is 
fiberwise invariant. We verify that it is of the form
\begin{align} \label{ctyri}
\U(u_0):=-{1 \over 2}F(p_0(u_0),u_0).
\end{align}
Actually, for pairs $(p_0(u_0),u_0) \in M \times \ba_0$ we have 
$\U(\va_0(p_0(u_0),u_0))=-\U(u_0)$ because in the fiber over its center,
the symmetry is only right multiplication by a suitable element acting 
as $-$id on $\fg_1$, and the fact follows from the equivariancy.

Now we show that under the condition (ii), function $\U$ solves the 
equation on $M \times \ba_0$ or, equivalently, the fiberwise invariant 
Weyl structure $\hat \si = \si \cdot \exp \U$ is invariant.
We substitute $(\ref{ctyri})$ into the equation $(\ref{tri})$.
The formula is then of the form
\begin{align} \label{pet}
F(x,u_0)= -{1 \over 2}F(p_0(\va_0(x,u_0)),\va_0(x,u_0)) +{1 \over 2}F(p_0(u_0),u_0)
\end{align}
and this has to hold for each $x \in M$ and $u_0 \in \ba_0$. 
If we suppose $u_0 \in p_0^{-1}(y)$ for an arbitrary fixed $y \in M$,
we simply want to verify that
\begin{align} \label{sest}
F(x,u_0)= -{1 \over 2}F(s_x(y),\va_0(x,u_0)) +{1 \over 2}F(y,u_0).
\end{align}
Here we use the equation $s_x \circ s_y=s_{s_x(y)} \circ s_x$. 
We apply (the coverings of) the left and right hand side of this equation 
on $\si$. The left hand side is of the form 
\begin{align*}
\va(x,\va(y,\si(u_0)))&= \va\big(x,\si(\va_0(y,u_0))\cdot \exp F(y,u_0)\big)=
\\&\va \big(x,\si(\va_0(y,u_0))\big) \cdot \exp F(y,u_0)=
\\&\si(\va_0(x,\va_0(y,u_0))) \cdot \exp F(x,\va_0(y,u_0))\cdot \exp F(y,u_0)
.\end{align*}
The right hand side is of the form 
\begin{align*}
\va(s_x(y),\va(x,\si(u_0)))&= \va\big(s_x(y),\si(\va_0(x,u_0))\cdot \exp F(x,u_0)\big)=
\\&\va\big(s_x(y),\si(\va_0(x,u_0))\big)\cdot \exp F(x,u_0)=
\\&\si(\va_0(s_x(y),\va_0(x,u_0))) \cdot \exp F(s_x(y),\va_0(x,u_0))\cdot \exp F(x,u_0)
.\end{align*}
Because we have an effective geometry, each morphism determines uniquely its covering. 
Thus the equality $s_x \circ s_y=s_{s_x(y)} \circ s_x$ implies the equality of the 
coverings 
\begin{align*}
\va(x,\va(y,\si(u_0)))&= \va(s_x(y),\va(x,\si(u_0))),\\
\va_0(x,\va_0(y,u_0))&=\va_0(s_x(y),\va_0(x,u_0))
.\end{align*}
Comparing of the $\exp$ parts then gives
$$
F(x,\va_0(y,u_0))+F(y,u_0)= F(s_x(y),\va_0(x,u_0))+F(x,u_0).
$$
Moreover, the first term on the left hand side can be rewritten as 
$ F(x,\va_0(y,u_0))
= -F(x,u_0)
$
because $u_0 \in p_0^{-1}(y)$ and $\va_0(y,u_0)$ is simply the covering of the symmetry 
at $y$ applied on a point from the fiber over $y$.
We get 
$$
-F(x,u_0)+F(y,u_0)= F(s_x(y),\va_0(x,u_0))+F(x,u_0)
$$
and some arrangements give
$$
{1 \over 2}F(y,u_0)-{1 \over 2}F(s_x(y),\va_0(x,u_0))= F(x,u_0)
,$$
which is exactly $(\ref{sest})$. 
Thus we have verified that $(\ref{ctyri})$ solves $(\ref{tri})$ and thus 
$\hat \si=\si \cdot \exp \U$ is the invariant Weyl structure.
The uniqueness of the fiberwise invariant Weyl structure implies that this 
does not depend on the choice of $\si$ we started with.
\end{proof}

\subsection{Invariant connections}\label{4.3}
Let us formulate now the main result on symmetric $|1|$--graded geometries:
\begin{thm*}
Let $(\ba \rightarrow M,\om)$ be a $|1|$--graded geometry and 
let $S: M \times M\rightarrow M$, $S(x,y)=s_x(y)$ be a system of symmetries
such that 
\begin{enumerate}
\item[(i)]  $S$ admits fiberwise invariant Weyl structure,
\item[(ii)] $s_x \circ s_y \circ s_x = s_{s_x(y)}$ holds for each $x,y \in M$. 
\end{enumerate}
Then there is an admissible affine connection $\na$ which is invariant with respect 
to all symmetries from the system $S$. In particular, $(M,\na)$ is an affine symmetric space. 
\end{thm*}

\begin{proof}
Everything follows from the Proposition \ref{4.2}. Under the assumptions, there is a 
Weyl structure $\hat \si$ invariant with respect to all symmetries. Then the corresponding
Weyl connection $\na$ is invariant with respect to all symmetries 
(and clearly respects the geometrical structure). 
Then the pair $(M,\na)$ is an affine symmetric space for the symmetries from the system $S$. 
\end{proof}
Let us point out here that for a system $S$, the invariant connection $\na$ is given uniquely. 
This follows directly from the uniqueness of the (fiberwise) invariant Weyl structure.
In fact, the question of the uniqueness of the mapping $\hat \si$ constructed in 
Lemma \ref{4.1} is independent from the question on its smoothness. Such a mapping 
$\hat \si: \ba_0 \rightarrow \ba$ is 
always given uniquely, but if it is not smooth, then it does not define a connection.

Let us also remark that the Theorem in particular implies that the system $S$ satisfying 
both conditions is smooth.

\subsection{Smooth systems and homogeneity}
Proposition \ref{4.2} and Theorem \ref{4.3} in particular apply in the case of a 
smooth system $S$ on a $|1|$--graded geometry $(\ba \rightarrow M, \om)$. 
Such a system  clearly admits fiberwise invariant Weyl structure and thus the 
condition (i) is automatically satisfied.
Then $M$ has to be homogeneous, see \ref{3.2}. 
We have the following trivial observation:
\begin{cor*}
If a $|1|$--graded geometry $(\ba\rightarrow M,\om)$ carries a smooth system of symmetries 
$S$ such that $s_x \circ s_y \circ s_x=s_{s_x(y)}$ holds for each $x,y \in M$, then 
$(M, S)$ form an affine symmetric space.
\end{cor*}
But Theorem \ref{4.3} says that weaker conditions are sufficient to get 
an affine symmetric space.  
Let $(\ba \rightarrow M,\om)$ be a $|1|$--graded geometry with 
a system $S: M \times M\rightarrow M$, $S(x,y)=s_x(y)$ as in Theorem \ref{4.3}, i.e.
there is a Weyl connection $\na$ such that $(M,\na)$ is an affine symmetric 
(for the system $S$). Then, in particular, $M$ is homogeneous because the group 
of affine transformations of $\na$ 
(which involves all symmetries from $S$) acts transitively on $M$. 
Because all symmetries are morphisms of the $|1|$--graded geometry, 
$(\ba \rightarrow M,\om)$ has to form a homogeneous $|1|$--graded geometry
in the sense of Section \ref{3}.

Let us finally remark that one can see e.g. \cite{ja-srni08} and references 
therein for an explicit description of such an example. We discuss there a Grassmannian 
symmetric geometry as a space of chains.

\subsection{Uniqueness of symmetries}
Proposition \ref{2.3} says that many of the interesting geometries
allow at most one symmetry at a point with a non--zero curvature. 
For these geometries, the condition (ii) in Theorem \ref{4.3} is trivially 
satisfied and we have the following consequence:
\begin{cor*}
Let $S: M \times M\rightarrow M$, $S(x,y)=s_x(y)$ be a system of symmetries on a
non--flat geometry of projective type, almost quaternionic type 
or conformal type with positive or negative definite signature 
such that $S$ admits fiberwise invariant Weyl structure.
Then there is an admissible affine connection $\na$ such that $(M,\na)$ 
is an affine symmetric space. 
In particular, homogeneous non--flat symmetric geometries of the latter types 
are affine symmetric spaces.
\end{cor*}
Except for some strange examples, the description of the latter types
reduces to the classical case. Otherwise, if the system of symmetries is sufficiently 
nice i.e. admits fiberwise invariant Weyl structure, than we get an affine 
symmetric space and the previous ideas says that there are no other more 
interesting examples (with smooth system of symmetries). Moreover, there are 
exceptions with non--smooth systems of symmetries as in \ref{3.3}. 
(See also \cite{P}).

\end{document}